\documentclass[11pt]{article}
\usepackage[linktocpage=true]{hyperref}
\usepackage{setspace}
\usepackage{amssymb, amsmath, amsthm, graphicx,mathrsfs}
\usepackage[left=1in,top=1in,right=1in,bottom=1in]{geometry}

\usepackage{lipsum}
\usepackage{titlesec}
\makeatletter
\def\thm@space@setup{
  \thm@preskip=0.4cm
  \thm@postskip=\thm@preskip}

\titleformat{\subsection}[runin]
{\normalfont\bfseries}{\thesubsection}{1em}{}

\titleformat{\subsubsection}[runin]
{\normalfont\bfseries}{\thesubsubsection}{0.5em}{}

\let\OLDthebibliography\thebibliography
\renewcommand\thebibliography[1]{
  \OLDthebibliography{#1}
  \setlength{\parskip}{0pt}
  \setlength{\itemsep}{0pt plus 0.3ex}
}
      \theoremstyle{plain}
      \newtheorem{theorem}{Theorem}%[section]

      \newtheorem{corollary}[theorem]{Corollary}

            \theoremstyle{definition}

      \theoremstyle{remark}

\def\o{\omega}
\begin{document}
\title{The number of trees in a graph}

\medskip

\author{
Dhruv Mubayi \thanks{Department of Mathematics, Statistics, and Computer Science, University of Illinois, Chicago, IL, 60607 USA.  Research partially supported by NSF grant DMS-1300138. Email: {\tt mubayi@uic.edu}}\and Jacques Verstra\" ete \thanks{
Department of Mathematics, University of California, San Diego, CA 92093-0112. Research supported by NSF
grant DMS-1362650. Email: {\tt jacques@ucsd.edu}} }

\maketitle

\begin{abstract}

Let $T$ be a tree with $t$ edges. We show
that the number  of isomorphic (labeled) copies of $T$ in a graph $G = (V,E)$
of minimum degree at least $t$ is at least
\[2|E| \prod_{v \in V}  (d(v) - t + 1)^{\frac{(t-1)d(v)}{2|E|}}.\] Consequently,  any $n$-vertex graph of average degree  $d$ and minimum degree at least $t$ contains at least
 $$nd(d-t+1)^{t-1}$$
 isomorphic (labeled) copies of $T$.
 This answers a question of~\cite{Detal} (where the above statement was proved when $T$ is the path with three edges) while extending an old result of Erd\H os and Simonovits~\cite{ES}.
\end{abstract}

\section{Introduction}

Let $T$ be a  $t$-edge tree and $G=(V,E)$ be a graph with minimum degree at least $t$. In this note we consider the question of how  many (isomorphic) copies of $T$ we can find in $G$. More precisely, if $V(T)=\{x_1, \ldots, x_{t+1}\}$, then we wish to count the number of injections $\phi: V(T) \rightarrow V$ such that $\phi(u)\phi(v)\in E$ for every edge $uv$ of $T$. This is a basic question in combinatorics, for example, the simple lower bound $\sum_{v \in V}t!{d(v) \choose t}$ in the case when $T$ is a star is the main inequality needed for a variety of fundamental problems in extremal graph theory.

 A natural way to count walks of length $t$ in a graph $G$ is to add up the entries
of $A^t$, where $A$ is the adjacency matrix of $G$. The Blakley-Roy~\cite{BR} inequality
uses linear algebra to show that the number of walks of length $t$ is at least
$nd^t$ in any graph of average degree $d$ with $n$ vertices (in fact the inequality
is a more general statement about inner products). Another approach to counting
walks, and more generally homomorphisms of trees, was used by Sidorenko, using an analytic method and the tensor power trick~\cite{sidorenko}.
Erd\H os and Simonovits~\cite{ES} proved that in a graph with average degree $d$, the number of walks of length $t$ that repeat a vertex is a negligible proportion of the total number of walks of length $t$ as $d \rightarrow \infty$. Consequently, their result implies that in a graph of average degree $d$ with $n$ vertices there are
at least $(1 - o(1))d^tn$ paths with $t$ edges as $d\rightarrow \infty$.
On the other hand, in~\cite{Detal} the following lower bound for the number of homomorphic copies of $T$ in a graph $G=(V,E)$ is proved, where a homomorphic copy is a (not necessarily injective) function $\phi: V(T) \rightarrow V$ such that $\phi(u)\phi(v)\in E$ for every edge $uv$ of $T$:
\begin{equation} \label{detal}
2|E|\prod_{v \in V} d(v)^{\frac{(t-1)d(v)}{2|E|}}.\end{equation}

Combining (\ref{detal}) with the result of~\cite{ES},
we obtain that the number of isomorphic copies of a path $T$ in $G$ is at least
$$(1-o(1))2|E|\prod_{v \in V} d(v)^{\frac{(t-1)d(v)}{2|E|}}.$$
The result of Erd\H os and Simonovits~\cite{ES} does not
give a precise expression for the $o(1)$ error term above
(although presumably it could be worked out from their proof).

In~\cite{Detal} the following more precise lower bound was given in the case when $T=P_3$, the path with three edges, and $G$ has minimum degree at least $3$:
\begin{equation} \label{p3}
2|E|\prod_{v \in V} (d(v)-2)^{\frac{2d(v)}{2|E|}}.\end{equation}
The authors in~\cite{Detal} asked whether a bound similar to (\ref{detal}) and (\ref{p3}) could be proved for the number of isomorphic copies of a tree $T$ in $G$ assuming that $G$ has sufficiently large minimum degree. The spirit of the question was to obtain a bound that is a convex function of the degrees of the vertices (and in particular whose unique minimum occurs when $G$ is regular). Here we provide such a bound that generalizes (\ref{p3}).

\begin{theorem} \label{loc}
Let $T$ be a tree with $t$ edges and $G$ be an $n$-vertex graph with  average degree $d$ and minimum degree at least $t$. Then the number of isomorphic (labeled) copies of
$T$ in $G$ is at least
$$nd \prod_{v \in V}  (d(v) - t + 1)^{\frac{(t-1)d(v)}{nd}}.$$
\end{theorem}
A consequence of this  is the following lower bound
in terms of the average degree in $G$.

\begin{corollary} \label{cor}
Let $T$ be a tree with $t$ edges and $G$ be an $n$-vertex graph with  average degree $d$ and minimum degree at least $t$. Then the number of isomorphic (labeled) copies of
$T$ in $G$ is at least
$$nd(d-t+1)^{t-1}.$$
\end{corollary}
 Indeed, Corollary~\ref{cor}  follows immediately from  Theorem~\ref{loc} by applying Jensen's inequality to the function $f(x)=(t-1)x\log (x - t + 1)$ which
is convex for $x \geq t$. Note also that the Corollary is nearly sharp as shown by complete graphs.  Indeed, if $G$ is the $n$ vertex graph of disjoint cliques, each with $d+1$ vertices and $d \ge t$, then the number of copies of $T$ in $G$ is $nd(d-1)\cdots (d-t+1)$.

The proof of Theorem~\ref{loc} uses the ideas first introduced by Alon, Hoory and Linial~\cite{AHL}, and subsequently developed in~\cite{Detal} to count homomorphisms of $T$ in $G$. 

\section{Proof of Theorem \ref{loc}}

We start with a graph $G$ of minimum degree at least $t$ and a tree $T$ with $t$ edges. Let $\Omega$ be the set of
all isomorphic copies of $T$ in $G$. In other words, $\Omega$ is the set of injections $\xi: V(T) \rightarrow V(G)$ such that $\xi(u)\xi(v)\in E(G)$ for every $uv \in E(T)$. Label the vertices of $T$ by first fixing a leaf $x_1$ and then labeling vertices $x_2, x_3, ..$ such that for any $j>1$ there is a unique $f(j)<j$ such that $x_jx_{f(j)} \in E(T)$. We could, for example, label the vertices using Breadth-First Search
or Depth-First Search. Let us call such a labeling of $T$ {\em good}.

\medskip

We consider {\it oriented} isomorphisms $\phi~:~V(T) \rightarrow V(G)$ which
can be constructed as follows.  Start with an arbitrary (directed) edge $v_1v_2 \in E(G)$ and map $x_1$ to $v_1$ and $x_2$ to $v_2$.
 Once $x_1,x_2,\dots,x_i \in V(T)$ are embedded as $\o_1,\o_2,\dots,\o_i \in V(G)$, i.e., $\phi(x_j)
= \omega_j$ for $j \leq i$,  then choose an arbitrary neighbor $\o_{i+1}$ of $\o_{f(i+1)}$ outside $\{\o_1,\o_2, \ldots, \o_i\}$ and embed $x_{i+1}$ as  $\o_{i+1}$. This gives us a natural probability on the sample space $\Omega$ of isomorphic copies of $T$ in $G$, with associated probability measure $\mathbb P$. For convenience, given $\omega \in \Omega$, we let $\omega_i$ denote the $i$th vertex of $\omega$ in the embedding.
This probability measure is defined on a specific isomorphic copy $\omega \in \Omega$ of $T$ in $G$ by
\[ \mathbb P(\omega) = \frac{1}{nd} \prod_{i = 2}^{t} \frac{1}{|N(\omega_{f(i+1)}) \backslash \{\omega_{1},\omega_{2},\dots,\omega_{i}\}|}.\]
Since $|N(\omega_{f(i+1)}) \backslash \{\omega_{1},\omega_{2},\dots,\omega_{i}\}| \geq d(\omega_{f(i + 1)}) - t + 1$,
\[ \mathbb P(\omega)  \leq \frac{1}{nd} \prod_{i = 2}^{t} \frac{1}{d(\omega_{f(i+1)}) - t + 1} := p(\omega).\]
Let $d$ be the average degree of $G$ and $n$ the number of vertices. Then, by the inequality of arithmetic and geometric means (using that $\sum_{\o}\mathbb P(\o)=1$),
\[ |\Omega| \geq \prod_{\omega \subset G} \mathbb P(\omega)^{-\mathbb P(\omega)} \geq \prod_{\omega \subset G} p(\omega)^{-p(\omega)}
= nd \prod_{\omega \subset G} \prod_{i = 2}^{t} (d(\omega_{f(i+1)}) - t + 1)^{p(\omega)}.\]
Interchanging the products we get
\[ |\Omega| \geq nd \prod_{i = 2}^{t} \prod_{\omega \subset G} (d(\omega_{f(i+1)}) - t + 1)^{p(\omega)}. \]
A term in the product above of the form $d(v)-t+1$ appears
when $v$ is the $i$th vertex of some $\o \in \Omega$, for some $i : 2\le i \le t$. Therefore, we have
\begin{equation}
\label{products}
 |\Omega|\ge nd \prod_{i = 2}^t \prod_{v \in V} (d(v) - t + 1)^{g_i(v)}\end{equation}
where
\[ g_i(v) := \sum_{{\omega \subset G}\atop{\omega_i = v}} p(\omega).\]
The key part of the proof is to show 
\begin{equation} \label{giv} 
g_i(v) \geq \frac{d(v)}{nd}.
\end{equation} We note that here is where our proof differs from the previous works~\cite{AHL, Detal}.
Those papers dealt with homomorphisms instead of isomorphisms, so there was no need to avoid  previously embedded  vertices of $T$, and the corresponding probability distribution in that setting is
\[ \mathbb P'(\omega) = \frac{1}{nd} \prod_{i = 2}^{t} \frac{1}{d(\omega_{f(i+1)})}.\]
 Moreover, if we use the probabilitiy measure $\mathbb P'$ (instead of the function $p$ which is not a probability measure), then (\ref{giv}) actually holds with equality  essentially because the Markov chain associated with the distribution $\mathbb P'$ is reversible.  This is not true in our case, and there are  constructions  showing that in our situation, $g_i(v)> d(v)/nd$ is possible. Consequently,  the argument showing (\ref{giv}) in our situation is more delicate. 

 To this end, we will prove (\ref{giv}) by proving the following stronger statement by induction on $t$: Given a $t$-edge tree $T$ with good labeling $x_1, \ldots, x_{t+1}$ and associated function $f$, an $n$-vertex graph $G=(V,E)$ with average degree $d$, and $1 \le i \le t+1$, we have
$g_i(v) \geq d(v)/nd$.  Note that we have included $i=1$ and $i=t+1$ in this statement as this will be needed in the induction argument that we will use.

The case $t=1$ is trivial (for both $i=1$ and $i=2$) so assume that $t>1$. Let us first assume that $i<t+1$. Let $T'=T-x_{t+1}$ be the tree obtained from $T$ by deleting the leaf $x_{t+1}$, let $\o^-=\o_1, \ldots, \o_t$ and $N=N(\o_{f(t+1)})\setminus\{\o_1, \ldots \o_t\}$ so that $|N|\ge d(\o_{f(t+1)})-t+1$. Then
\begin{eqnarray*}
g_i(v) &=& \sum_{\o^-: \o_i=v}\sum_{\o_{t+1} \in N}p(\o) \\
&\ge& \sum_{\o^-: \o_i=v} (d(\o_{f(t+1)})-t+1)p(\o) \\
&\ge& \sum_{\o^-: \o_i=v}\frac{1}{nd}\prod_{i=2}^t \frac{1}{d(\o_{f(i+1)})-t+2}.
\end{eqnarray*}

Finally, we note that the rightmost expression is precisely $g_i(v)$ for the tree $T'$ which has $t-1$ edges. So by induction it is at least $d(v)/nd$ as required.

We now suppose that $i=t+1$. Given a copy $\o=\o_1,\ldots, \o_{t+1}$ of $T$ in $G$ with $\phi(x_i)=\o_i$ as usual, let us relabel the vertices with $z=z_1, \ldots,  z_{t+1}$ such that $z_1=\o_{t+1}, z_{t+1}=\o_1$, and $z$ is a good labeling of $\o=\phi(T)$. Note that this is clearly possible as we may just produce a good labeling of $\o$ beginning with $\o_{t+1}$ and ending with $\o_1$ (recall that $x_1$ is a leaf of $T$). As before, define
$$p(z):=\frac{1}{nd} \prod_{j = 2}^{t} \frac{1}{d(z_{f(j+1)}) - t + 1}.$$
Now we make the crucial observation that $p(\o)=p(z)$.
To see this, observe that
$$p(\o)=\frac{1}{nd} \prod_{j = 2}^{t} \frac{1}{d(\omega_{f(j+1)}) - t + 1}=\frac{1}{nd} \prod_{j = 2}^{t} \left(\frac{1}{d(\o_{j}) - t + 1}\right)^{d_T(x_j)-1}.$$
as each term is counted once for each child of  the corresponding vertex as the good labeling is constructed.
As $\{z_2, \ldots, z_t\}=\{\o_2, \ldots, \o_t\}$, we
also obtain
$$p(\o)=\frac{1}{nd} \prod_{j = 2}^{t} \left(\frac{1}{d(\o_{j}) - t + 1}\right)^{d_T(x_j)-1}=\frac{1}{nd} \prod_{j =2}^{t} \left(\frac{1}{d(z_{j}) - t + 1}\right)^{d_T(x_{\pi(j)})-1}=p(z),$$
where $\pi$ is the permutation on $t-1$ elements such that $z_j=x_{\pi(j)}$ for all $2 \le j \le t$. Consequently,
$$g_{t+1}(v)=\sum_{\omega: \omega_{t+1} = v} p(\omega)
=\sum_{z: z_1=v}p(\o)=
\sum_{z: z_1=v}p(z) =g_1(v) \ge \frac{d(v)}{nd}.$$

Inserting this into  (\ref{products}) we get
\[ |\Omega| \geq nd \prod_{v \in V} (d(v) - t + 1)^{\frac{(t - 1)d(v)}{nd}}.\]
This proves the theorem. \qed

\section{Concluding Remarks}

$\bullet$ If the maximum degree of the subgraph induced by any tree with $t$ edges is $k$, then the above proof gives a better bound:

\begin{corollary}
Fix a tree $T$ with $t$ edges. Let $G = (V,E)$ be an $n$-vertex graph such that copy of $T$ in $G$ induces a subgraph of maximum degree at most $k$,
and such that $G$ has minimum degree at least $k$. Then the number of isomorphic copies of $T$  in $G$ is at least
\[ 2|E| \prod_{v \in V} (d(v) - k + 1)^{\frac{(t - 1)d(v)}{2|E|}}.\]
\end{corollary}

$\bullet$ We were not able to decide if the following statement is true (even for paths):

\begin{center}
\parbox{5.7in}{
{\sl Fix a tree $T$ with $t$ edges. The number of isomorphic labeled copies of $T$  in an $n$-vertex graph of large enough minimum degree and average degree $d$ is at least $nd(d-1)\cdots (d-t+1)$.}}
\end{center}

This statement if true would be best possible, since a graph consisting of disjoint cliques of order $d + 1$ has
average degree $d$ and exactly $nd(d - 1)\cdots (d - t + 1)$ isomorphic copies of any tree with $t$ edges.

\end{document}